# The block spectrum of RNA pseudoknot structures

**Thomas J. X. Li** · **Christina S. Burris** ·
**Christian M. Reidys**

**Abstract**  In this paper we analyze the length-spectrum of blocks in $\gamma$-structures. $\gamma$-structures are a class of RNA pseudoknot structures that plays a key role in the context of polynomial time RNA folding. A $\gamma$-structure is constructed by nesting and concatenating specific building components having topological genus at most $\gamma$. A block is a substructure enclosed by crossing maximal arcs with respect to the partial order induced by nesting. We show that, in uniformly generated $\gamma$-structures, there is a significant gap in this length-spectrum, i.e., there asymptotically almost surely exists a unique longest block of length at least $n - O(n^{1/2})$ and that with high probability any other block has finite length. For fixed $\gamma$, we prove that the length of the longest block converges to a discrete limit law, and that the distribution of short blocks of given length tends to a negative binomial distribution in the limit of long sequences. We refine this analysis to the length spectrum of blocks of specific pseudoknot types, such as H-type and kissing hairpins. Our results generalize the rainbow spectrum on secondary structures by the first and third authors and are being put into context with the structural prediction of long non-coding RNAs.

Thomas J. X. Li
Biocomplexity Institute of Virginia Tech
Blacksburg, VA 24061, USA
E-mail: thomasli@vt.edu

Christina S. Burris
Biocomplexity Institute of Virginia Tech
Blacksburg, VA 24061, USA
E-mail: christie.burris@gmail.com

Christian M. Reidys
Biocomplexity Institute of Virginia Tech
Blacksburg, VA 24061, USA
E-mail: duck@santafe.edu

# 1 Introduction

Ribonucleic acid (RNA) plays an important role in various biological processes within cells, ranging from catalytic activity to gene expression. High throughput sequencing technique has revealed a large number of non-coding RNA transcripts whose functions in biological processes are beginning to be explored. In particular, long non-coding RNAs (lncRNAs) are emerging as an integral functional component of the human transcriptome and have attracted substantial attention in the past few years (Iyer et al, 2015).

An RNA molecule folds into a helical configuration of its primary sequence by forming hydrogen bonds between pairs of nucleotides according to Watson-Crick and wobble base-pairing rules. These structures are often key to understanding their functions within cells such as: transcription and translation (McCarthy and Holland, 1965), catalyzing reactions (Kruger et al, 1982), gene regulation (Eddy, 2001).

The most prominent class of coarse grained RNA structures are the RNA secondary structures. They encode the bonding information of the nucleotides irrespective of the actual spacial embedding. More than three decades ago, Waterman and his coworkers pioneered the combinatorics and prediction of RNA secondary structures (Waterman, 1978, 1979; Smith and Waterman, 1978; Howell et al, 1980; Schmitt and Waterman, 1994; Penner and Waterman, 1993). Represented as a *diagram* by drawing its sequence on a horizontal line and each base pair as an arc in the upper half-plane, RNA secondary structure contains no crossing arcs (wo arcs $(i_1, j_1)$ and $(i_2, j_2)$ cross if the nucleotides appear in the order $i_1 < i_2 < j_1 < j_2$ in the primary sequence).

In fact, it is well-known that there exist cross-serial interactions, called pseudoknots in RNA (Westhof and Jaeger, 1992), see Fig. 1. RNA structures with cross-serial interactions are of biological significance, occur often in practice and are found to be functionally important in tRNAs, RNAseP (Loria and Pan, 1996), telomerase RNAs (Staple and Butcher, 2005; Chen et al, 2000), and ribosomal RNAs (Konings and Gutell, 1995). Cross-serial interactions also appear in plant viral RNAs and *in vitro* RNA evolution experiments have produced pseudoknotted RNA families, when binding HIV-1 reverse transcriptase (Tuerk et al, 1992).

The key to organize and filter structures with cross-serial interactions is to introduce topology. The idea is simple: instead of drawing a structure in the plane (sphere) we draw it on more sophisticated orientable surfaces. The advantage of this is that this presentation allows to eliminate any cross-serial interactions. The topology of RNA structures has first been studied in Penner and Waterman (1993); Penner (2004) and the classification of RNA structures including pseudoknots in terms of the topological genus of an associated fat-graph via matrix theory in Orland and Zee (2002); Vernizzi et al (2005); Bon et al (2008). Andersen et al (2013) study topological RNA structures of higher genus and associate them with Riemann's Moduli space in Penner (2004).

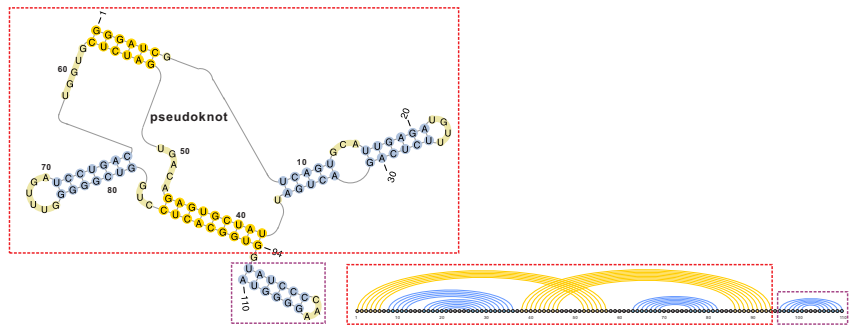

**Fig. 1** The pseudoknot structure and diagram representation of ribox02 (Tsukiji et al, 2003). It contains one block of H-type having length 94, and one block of T-type having length 15 (highlighted by dotted boxes).

The topology of RNA pseudoknot structures has been further translated into an efficient dynamic programming algorithm (Reidys et al, 2011). This algorithm *a priori* folds into a novel class of pseudoknot structures, the $\gamma$-structures. $\gamma$-structures differ from pseudoknotted RNA structures of fixed topological genus (Orland and Zee, 2002; Bon et al, 2008). They are composed by irreducible components whose individual genus is bounded by $\gamma$ and contain no bonds of length one (1-arcs), see Section 2 for details. RNA pseudoknot structures are thus filtered by just one parameter $\gamma$. RNA secondary structure fits seamlessly into this classification, since these are exactly 0-structures, i.e., structures consisting of components of genus zero (noncrossing arcs). Han et al (2014) study the combinatorics of $\gamma$-structures, by deriving simple asymptotic formulas for the number of $\gamma$-structures. Li and Reidys (2013) show a central limit theorem for the distribution of topological genus in $\gamma$-structures of given length.

This paper is motivated by the recent work on the rainbow spectrum of RNA secondary structures (Li and Reidys, 2018). A rainbow in a secondary structure is a maximal arc with respect to the partial order induced by nesting, i.e. the closing arc of a stem-loop. The length of a rainbow $(i, j)$, defined as $j - i$, reflects the size of the corresponding stem-loop. Li and Reidys (2018) show that, in uniformly generated RNA secondary structures, there exists a unique longest rainbow of expected length $n - O(n^{1/2})$. For pseudoknot structures, arcs are allowed to cross and become associated with each other via a sequence of pairwise crossing arcs. These associated arcs form larger substructures, see Fig. 1. A particular class of such substructures are blocks. A block is a substructure enclosed by associated maximal arcs. In this paper we study the length spectrum of blocks (block spectrum) in $\gamma$-structures, generalizing the results of Li and Reidys (2018) on the block spectrum in 0-structures.

We further classify each block into different types based on the crossing pattern of its enclosing maximal arcs. While a block is called T-type if it is enclosed by a rainbow that does not cross with any other arcs, it also admits different crossing patterns, such as H-type, kissing hairpin (K-type), 3-knot (L-type) and 4-knot pseudoknots (M-type), see Fig. 1. We shall enrich our results on the block spectrum by considering blocks of different types.

The key results of this paper are the following:

1. in uniformly generated $\gamma$-structures the length of the longest block tends, in the limit of long sequences, to a discrete limit law, having an expectation value $n - O(n^{1/2})$. That is, there exists a unique longest block,
2. the probability of the unique longest block to be of certain type has a nonzero limit, which depends on $\gamma$ and its type,
3. with high probability any other block has finite length, $k$,
4. in the limit of long sequences, the distribution of blocks of any type having length $k$ tends to a negative binomial distribution,
5. all the above results apply to structures with given minimum stack- and arc-length constraints.

Our results are important in the context of the structural prediction of lncRNAs. Due to the high computational cost, many prediction algorithms for the lncRNA structure adopt a sliding window technique, i.e., restricting the length of the blocks by a predefined constant (Reeder et al, 2007). Our results show that, however, sliding window methods are incompatible with the unique giant block that we observe in RNA structures.

Our analysis is based on two concepts. The first is a new decomposition for $\gamma$-structures into blocks, generalizing the standard decomposition of secondary structures of Waterman (1978). This decomposition allows us to derive generating functions for $\gamma$-structures, either with a restricted length of their longest block or with a restricted type of block. These generating functions lay the basis for the analysis of their corresponding limit distributions, either for the unique longest block or for the short blocks.

The second is the fact that different block types are characterized by irreducible shadows. Here, irreducible means that the diagram cannot be decomposed via nesting and concatenating, and a shadow is a diagram without unpaired vertices, arcs of length 1 and parallel arcs. Intuitively, an irreducible shadow can be viewed as the minimal building element of a $\gamma$-structure, and has been studied in Han et al (2014); Li and Reidys (2013). Augmenting the generating polynomial for irreducible shadows computed in Han et al (2014), we are able to generalize the block spectra analysis to blocks of fixed type.

This paper is organized as follows: In Section 2, we provide some basic facts of $\gamma$-structures. We provide in Section 3 a novel decomposition of $\gamma$-structures into blocks. As a result, we derive a new generating function for $\gamma$-structures and then extract its coefficient asymptotics. In Section 4, we compute the expectation and variance of the longest block in $\gamma$-structures and prove the discrete limit law. We then observe in Section 5 that with high probability we can restrict our analysis to blocks of finite length and proceed computing the associated limit distribution. In Section 6 we extend our analysis to different block types and we integrate our results in Section 7.

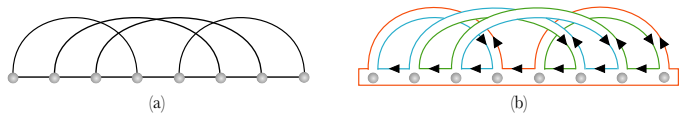

**Fig. 2** A diagram, its corresponding fatgraph and boundary components. (a) A diagram $G$ of genus 1. (b) Its corresponding fatgraph $\mathbb{G}$ represented by an orientable surface $F(\mathbb{G})$ with three boundary components (orange, blue, green).

## 2 Basic facts

An RNA secondary structure can be represented as a *diagram*, a labeled graph over the vertex set $\{1, \ldots, n\}$ whose vertices are arranged in a horizontal line and arcs are drawn in the upper half-plane. Clearly, vertices correspond to nucleotides in the primary sequence and arcs correspond to the Watson-Crick **A-U**, **C-G** and wobble **U-G** base pairs. The *length* of the structure is defined as the number of nucleotides. The length of an arc $(i, j)$ is defined as $j - i$. The backbone of a diagram is the sequence of consecutive integers $(1, \ldots, n)$ together with the edges $\{\{i, i+1\} \mid 1 \leq i \leq n-1\}$. We shall distinguish the backbone edge $\{i, i+1\}$ representing a phosphodiester bond, from the arc $(i, i+1)$, which we refer to as a 1-*arc*. Two arcs $(i_1, j_1)$ and $(i_2, j_2)$ are *crossing* if $i_1 < i_2 < j_1 < j_2$.

An RNA *secondary structure* is a diagram without 1-arcs and crossing arcs (Waterman, 1978). A *stack* of length $r$ is a maximal sequence of "parallel" arcs, $((i, j), (i+1, j-1), \ldots, (i+(r-1), j-(r-1)))$. A structure is *r-canonical* if it has minimum stack-length $r$.

We shall consider diagrams as fatgraphs, i.e., graphs together with a collection of cyclic orderings on the half-edges incident to each vertex. A fatgraph is obtained from a diagram by expanding each vertex to a disk and fattening the edges into (untwisted) ribbons such that the ribbons connect the disks in the given cyclic orderings. The specific drawing of a diagram $G$ with its arcs in the upper half-plane determines a collection of cyclic orderings on the half-edges of the underlying graph incident on each vertex, thus defining a corresponding fatgraph $\mathbb{G}$, see Fig. 2. Accordingly, each fatgraph $\mathbb{G}$ determines an associated orientable surface $F(\mathbb{G})$ with boundary (Loebl and Moffatt, 2008; Penner et al, 2010), which contains $G$ as a deformation retract (Massey, 1967), see Fig. 2. Fatgraphs were first applied to RNA secondary structures in Penner and Waterman (1993) and Penner (2004).

A diagram $\mathbb{G}$ hence determines a unique surface $F(\mathbb{G})$ (with boundary). Filling the boundary components with discs we can pass from $F(\mathbb{G})$ to a surface without boundary. The Euler characteristic, $\chi$, and genus, $g$, of this surface are given by $\chi = v - e + r$ and $g = 1 - \frac{1}{2}\chi$, respectively, where $v, e, r$ is the number of discs, ribbons and boundary components in $\mathbb{G}$ (Massey, 1967). The *genus* of a diagram is that of its associated surface without boundary.

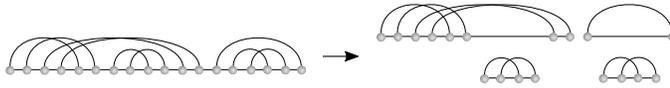

**Fig. 3** A 2-structure and its decomposition into components having genus at most 2.

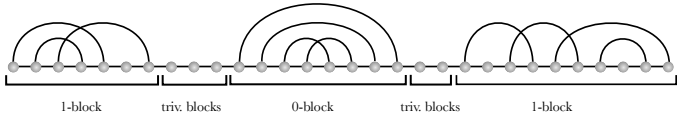

**Fig. 4** The decomposition of a $\gamma$-structure into its blocks: trivial blocks, 0-blocks and $\gamma$-blocks.

Two $G$-arcs $\alpha_1$ and $\alpha_k$ are *associated* if there exists a sequence of $G$-arcs $(\alpha_1, \alpha_2, \ldots, \alpha_k)$ such that $(\alpha_i, \alpha_{i+1})$ are crossing. Association is an equivalence relation and partitions the set of $G$-arcs into equivalence classes. For each equivalence class, $A$, a *component* is the induced diagram obtained from $G$ by removing all arcs not in $A$ and all isolated vertices. Therefore, any diagram $G$ without isolated vertices decomposes into a set of its components, see Fig. 3. Furthermore, the genus is additive with respect to this decomposition, i.e., the genus of the diagram is the sum of the genera of its components.

A $\gamma$-*structure* is a diagram $G$ without 1-arcs such that any $G$-component has genus at most $\gamma$, see Fig. 3. In particular, a 0-structure is an RNA secondary structure.

We define the partial order on the set of $G$-arcs to be $(i, j) \leq (i', j') \iff i' \leq i \ \wedge j \leq j'$. A *maximal component* is a component containing some maximal arc. An *exterior vertex* is an unpaired vertex $k$ such that there is no arc $(i, j)$ such that $i < k < j$. Clearly, the left- and rightmost endpoints of maximal components together with exterior vertices induce a partition of the backbone into subsequent intervals. A *block* is the induced diagram of $G$ over each such interval. By construction, a block is either trivially an exterior vertex or bounded by a maximal component. Thus, any diagram $G$ decomposes uniquely into a set of blocks, see Fig. 4. Again, the genus of the diagram is the sum of the genera of its blocks.

Each nontrivial block is characterized by its unique maximal component. A block is called a 0-*block* if its maximal component is of genus zero, i.e., a single arc. Otherwise, a block is a $\gamma$-*block* if its maximal component is of genus $1 \leq g \leq \gamma$. It turns out that the maximal component can be further reduced to an irreducible shadow defined as follows.

A *shadow* is a diagram without noncrossing arcs, isolated vertices and stacks of length greater than one. The shadow of a diagram is obtained by removing all noncrossing arcs, deleting all isolated vertices and collapsing all induced stacks to single arcs, see Fig. 5. Furthermore, projecting into the shadow does not affect genus. A diagram $G$ is *irreducible* if any two $G$-arcs are associated. Irreducibility is equivalent to the concept of primitivity introduced by Bon et al (2008). According to Andersen et al (2012), any irreducible shadow of genus $g$ has $2g \leq \ell \leq (6g - 2)$ arcs, i.e., there exist only finitely many irreducible shadows.

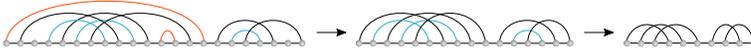

**Fig. 5** Shadows: the shadow is obtained by removing all noncrossing arcs (orange) and isolated points and collapsing all stacks (blue) and resulting stacks into single arcs.

## 3 $\gamma$-structures

We consider $\gamma$-structures subject to *minimum arc-length* and *minimum stack-length* constraint. This constraint is motivated by the fact that RNA structures having a minimum arc-length of four and a minimum stack length two or three are energetically more favorable. The former is a consequence of the rigidity of the molecules backbone (Stein and Waterman, 1979) and the latter reflects the fact that the main contribution of lowering free energy stems from electron delocalization between stacked bonds (Hunter and Sanders, 1990; Šponer et al, 2001, 2013).

Let $g_{\gamma,\lambda}^{[r]}(n)$ denote the number of $r$-canonical $\gamma$-structures over $n$ nucleotides with minimum arc-length $\lambda$. We shall simplify notation by writing $g_\gamma(n)$ instead of $g_{\gamma,\lambda}^{[r]}(n)$. Let $f_\gamma(n)$, $b_0(n)$ and $b_\gamma(n)$ denote the numbers of $r$-canonical blocks, 0-blocks and $\gamma$-blocks over $n$ nucleotides with minimum arc-length $\lambda$, respectively. The generating functions, filtered by sequence length, are given by

$$\mathbf{G}_\gamma(z) = \sum_{n \geq 0} g_\gamma(n) z^n, \quad \mathbf{F}_\gamma(z) = \sum_{n \geq 0} f_\gamma(n) z^n$$

$$\mathbf{B}_0(z) = \sum_{n \geq 1} b_0(n) z^n, \quad \mathbf{B}_\gamma(z) = \sum_{n \geq 1} b_\gamma(n) z^n.$$

Let $i_g(n)$ denote the number of irreducible shadows of genus $g$ with $n$ arcs, and define the generating polynomial $\mathbf{I}_g(z) = \sum_{n=2g}^{6g-2} i_g(n) z^n$. For instance for genus 1 and 2 we have

$$\mathbf{I}_1(z) = z^2 \left(1 + z\right)^2,$$
$$\mathbf{I}_2(z) = z^4 \left(1 + z\right)^4 \left(17 + 92\,z + 96\,z^2\right).$$

Han et al (2014) provides a recursion of $\mathbf{I}_g(z)$ for any given $g$. The recursion involves the generating function for unicellular fatgraphs of fixed genus (Harer and Zagier, 1986), whose computation is further facilitated by a recursion of certain coefficients $\kappa_g(n)$ (Li, 2014; Huang and Reidys, 2015).

The generating function $\mathbf{G}_\gamma(z)$ has been computed in Han et al (2014); Li and Reidys (2013) via an indirect inflation process from irreducible shadows to $\gamma$-structures using so-called "$\gamma$-matchings" and "$\gamma$-shapes" as intermediate objects.

Here we provide a novel approach to deriving a new and simple formula for $\mathbf{G}_\gamma(z)$. We utilize the decomposition of a $\gamma$-structure into a sequence of blocks and directly relate the construction of blocks to irreducible shadows. By doing so, we derive a functional equation which completely determines $\mathbf{G}_\gamma(z)$ and facilitates the asymptotic analysis for its coefficients $g_\gamma(n)$. As a byproduct, we obtain the generating functions for blocks, which are utilized in Section 6.

**Theorem 1** *For $\gamma, \lambda, r \geq 1$ with $\lambda \leq r+1$, the generating functions $\mathbf{G}_\gamma(z)$, $\mathbf{F}_\gamma(z)$, $\mathbf{B}_0(z)$ and $\mathbf{B}_\gamma(z)$ satisfy the functional equations*

$$\mathbf{G}_\gamma(z) = \frac{1}{1 - \mathbf{F}_\gamma(z)}, \tag{1}$$

$$\mathbf{F}_\gamma(z) = z + \mathbf{B}_0(z) + \mathbf{B}_\gamma(z), \tag{2}$$

$$\mathbf{B}_0(z) = \frac{z^{2r}}{1 - z^2}\Big(\mathbf{G}_\gamma(z) - \mathbf{B}_0(z) - \sum_{i=0}^{\lambda-2} z^i\Big), \tag{3}$$

$$\mathbf{B}_\gamma(z) = \mathbf{G}_\gamma(z)^{-1} \sum_{g \leq \gamma} \mathbf{I}_g\Big(\frac{z^{2r}\mathbf{G}_\gamma(z)^2}{1 - z^2 + z^{2r} - z^{2r}\mathbf{G}_\gamma(z)^2}\Big). \tag{4}$$

*In particular, there exists a polynomial $Q_\gamma(z, X)$ in $X$ of degree $(12\gamma - 2)$, such that*

$$Q_\gamma(z, \mathbf{G}_\gamma(z)) = 0. \tag{5}$$

The key idea to prove the system of functional equations of Theorem 1 is the following: any $\gamma$-structure can be decomposed into a sequence of blocks, and any block is either a single vertex, or a 0-block, or a $\gamma$-block. While a 0-block is formed of the stack containing the rainbow together with the enclosed $\gamma$-structure, a $\gamma$-block is characterized by the maximal component, see Fig. 6. The proof is presented in Section 8.

The implicit functional equation (5) completely characterizes the generating function $\mathbf{G}_\gamma(z)$. Using the singular implicit-function schema (see Flajolet and Sedgewick (2009)), we compute the singular expansion of $\mathbf{G}_\gamma(z)$ and the asymptotics of its coefficients, which are consistent with those in Han et al (2014); Li and Reidys (2013).

**Theorem 2 (Han et al (2014); Li and Reidys (2013))** *Suppose $1 \leq \gamma, \lambda, r \leq 4$ with $\lambda \leq r+1$. The dominant singularity $\rho_\gamma$ of $\mathbf{G}_\gamma(z)$ is the minimal positive, real solution of the resultant polynomial*

$$\Delta(z) = \mathbf{R}\left(Q_\gamma(z, X), \frac{\partial}{\partial X}Q_\gamma(z, X), X\right).$$

*The singular expansion of $\mathbf{G}_\gamma(z)$ is given by*

$$\mathbf{G}_\gamma(z) = \tau + \delta(\rho_\gamma - z)^{\frac{1}{2}} + \theta(\rho_\gamma - z) + O\big((\rho_\gamma - z)^{\frac{3}{2}}\big), \qquad as\ z \to \rho_\gamma,$$

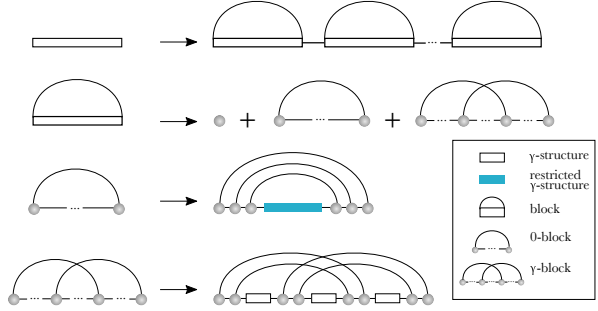

**Fig. 6** The decomposition of a $\gamma$-structure.

where $\tau = \mathbf{G}_\gamma(\rho_\gamma)$, $\delta$ and $\theta$ are constants, that can be explicitly computed. Furthermore, the coefficients of $\mathbf{G}_\gamma(z)$ satisfy

$$[z^n]\mathbf{G}_\gamma(z) = c\,n^{-\frac{3}{2}}\rho_\gamma^{-n}\left(1 + O(n^{-1})\right), \qquad as\ n \to \infty,$$

where $c$ is the positive constant $c = \delta\rho_\gamma^{\frac{1}{2}}\,\Gamma(-\frac{1}{2})^{-1}$.

In view of eq. (1), one can derive the asymptotics for blocks from the asymptotics for $\gamma$-structures.

**Corollary 1** Suppose $1 \le \gamma, \lambda, r \le 4$ with $\lambda \le r+1$. The dominant singularity of $\mathbf{F}_\gamma(z)$ is the same as that of $\mathbf{G}_\gamma(z)$, $\rho_\gamma$. The singular expansion of $\mathbf{F}_\gamma(z)$ is given by

$$\mathbf{F}_\gamma(z) = \tau' + \delta'\left(\rho_\gamma - z\right)^{\frac{1}{2}} + \theta'\left(\rho_\gamma - z\right) + O\left((\rho_\gamma - z)^{\frac{3}{2}}\right), \qquad as\ z \to \rho_\gamma,$$

where $\tau'$, $\delta' = \tau^{-2}\delta$ and $\theta'$ are constants. Furthermore, the coefficients of $\mathbf{F}_\gamma(z)$ satisfy

$$[z^n]\mathbf{F}_\gamma(z) = c'\,n^{-\frac{3}{2}}\rho_\gamma^{-n}\left(1 + O(n^{-1})\right), \qquad as\ n \to \infty, \tag{6}$$

where $c'$ is the positive constant $c' = \tau^{-2}c$.

## 4 The longest block

In this section, we shall show that the length of the longest block of uniformly generated $\gamma$-structures tends to a discrete limit law, having expected value $n - O(n^{1/2})$. These results are generalization of results on secondary structures in view of the fact that secondary structures are 0-structures  (Li and Reidys, 2018).

To this end, we analyze the random variable, $\mathbb{B}_{\gamma,n}$, representing the length of the longest block in a $\gamma$-structure of $n$ nucleotides. Clearly, $\mathbb{B}_{\gamma,n}$ is considered in the discrete probability space over all $\gamma$-structures of $n$ nucleotides with the uniform distribution $\mathbb{P}$, i.e., the distribution in which each structure $S$ has probability $\mathbb{P}(S) = \frac{1}{g_\gamma(n)}$.

In the following we derive first- and second-order information about $\mathbb{B}_{\gamma,n}$. This will later allow us to apply a large deviation result and imply that the random variable $n - \mathbb{B}_{\gamma,n}$ asymptotically almost surely (a.a.s.) converges to a discrete limit law.

**Lemma 1** *Suppose* $1 \leq \gamma, \lambda, r \leq 4$ *with* $\lambda \leq r + 1$. *The expectation and variance of* $\mathbb{B}_{\gamma,n}$ *are given by*

$$\mathbb{E}[\mathbb{B}_{\gamma,n}] = n - \alpha\, n^{\frac{1}{2}}\big(1 + o(1)\big), \qquad \mathbb{V}[\mathbb{B}_{\gamma,n}] = \beta n^{\frac{3}{2}}\big(1 + o(1)\big), \qquad as\ n \to \infty,$$

*where* $\alpha = 4\, c\, \tau^{-1}$ *and* $\beta = (1 - \frac{\pi}{4})\alpha$ *are positive constants, see Table 1.*

**Table 1** The values of $\alpha$ for $\gamma$-structures for $0 \leq \gamma \leq 3$ with minimum arc-length and minimum stack-length constraint $1 \leq \lambda, r \leq 4$ and $\lambda \leq r + 1$. The $\alpha$-values for 0-structures are computed in Li and Reidys (2018).

| | $\gamma = 0$ | | | | $\gamma = 1$ | | |
|---|---|---|---|---|---|---|---|
| | $r = 1$ | $r = 2$ | $r = 3$ | | $r = 1$ | $r = 2$ | $r = 3$ |
| $\lambda = 1$ | 1.954 | 2.804 | 3.431 | $\lambda = 1$ | 0.868 | 1.271 | 1.566 |
| $\lambda = 2$ | 1.687 | 2.483 | 3.096 | $\lambda = 2$ | 0.804 | 1.196 | 1.488 |
| $\lambda = 3$ | | 2.201 | 2.797 | $\lambda = 3$ | | 1.149 | 1.434 |
| $\lambda = 4$ | | | 2.529 | $\lambda = 4$ | | | 1.399 |

| | $\gamma = 2$ | | | | $\gamma = 3$ | | |
|---|---|---|---|---|---|---|---|
| | $r = 1$ | $r = 2$ | $r = 3$ | | $r = 1$ | $r = 2$ | $r = 3$ |
| $\lambda = 1$ | 0.640 | 0.941 | 1.162 | $\lambda = 1$ | 0.520 | 0.766 | 0.947 |
| $\lambda = 2$ | 0.601 | 0.896 | 1.115 | $\lambda = 2$ | 0.492 | 0.734 | 0.914 |
| $\lambda = 3$ | | 0.871 | 1.085 | $\lambda = 3$ | | 0.717 | 0.893 |
| $\lambda = 4$ | | | 1.066 | $\lambda = 4$ | | | 0.881 |

Our proof is analogous to that of Lemma 1 in Li and Reidys (2018), and outlined in Section 8.

**Remark:** Lemma 1 shows that the length of the longest block is $n - O(n^{\frac{1}{2}})$ with a standard deviation of $O(n^{\frac{3}{4}})$. As a result, the distribution of $\mathbb{B}_{\gamma,n}$ becomes for larger and larger $n$ more and more concentrated.

Table 1 shows that the parameter $\alpha$ decreases, if the structural complexity $\gamma$ or the minimum arc-length $\lambda$ increase, or the minimum stack-length $r$ decreases. Furthermore, the impact of $\lambda$ on $\alpha$-values is lower compared with $\gamma$ and $r$.

In Fig. 7, we contrast our asymptotic estimate of the expectation and the average length of the longest block in 1-structures.

**Theorem 3** *We have for any* $t > \frac{3}{4}$,

$$\lim_{n \to \infty} \mathbb{P}(n - \mathbb{B}_{\gamma,n} \geq \Omega(n^t)) = 0 \tag{7}$$

*and for any* $k = o(n)$

$$\lim_{n \to \infty} \mathbb{P}(n - \mathbb{B}_{\gamma,n} = k) = \tau^{-2}\, b_k\, \rho^k, \tag{8}$$

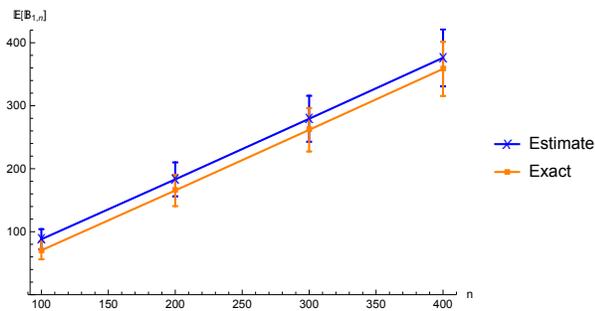

**Fig. 7** The longest block in 1-structures: we compare the expectation value and standard deviation (crosses) with the average length (squares, shifted by subtracting 20) computed by exact enumeration in all 1-structures. Minimum arc- and stack-length are $r = \lambda = 2$, $n = 100 \times i$. where $1 \leq i \leq 4$. Error bars represent one standard deviation.

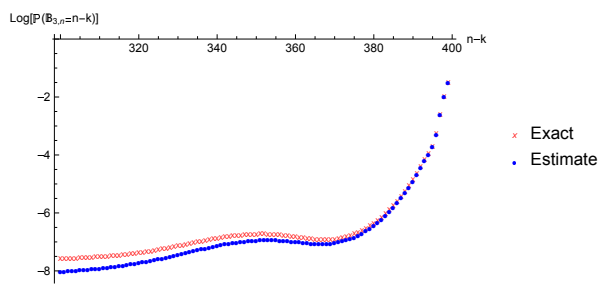

**Fig. 8** The longest block: we contrast the limit distribution (dots) with the distribution (crosses) of 3-structures. Minimum arc- and stack-length are $r = \lambda = 2$, sequence length is $n = 400$.

where $b_k = [z^k]\Phi'(\mathbf{F}_\gamma(z))$ and $\Phi(z) = \frac{1}{1-z}$. *Consequently the distribution of* $n - \mathbb{B}_{\gamma,n}$ *a.a.s. converges to a discrete limit law.*

Theorem 3 is a consequence of Lemma 1. The proof is a direct generalization of that of Theorem 3 in Li and Reidys (2018) to $\gamma$-structures, and is presented in Section 8.

In Fig. 8, we compare our theoretical result with the distribution of the length of the longest block computed by exact enumeration of all 3-structures.

Fig. 9 shows that the decrease of $\mathbb{P}(\mathbb{B}_{\gamma,n} = n-k)$, for increasing $k$, depends on $\gamma$, minimum stack- and arc-length $r, \lambda$. For $k > 10$, the probability of the longest block having length $n - k$ becomes smaller, when $\gamma$ or $\lambda$ increase. In contrast, we observe a positive correlation between $\mathbb{P}(\mathbb{B}_{\gamma,n} = n-k)$ and $r$, see Fig. 9 (Middle).

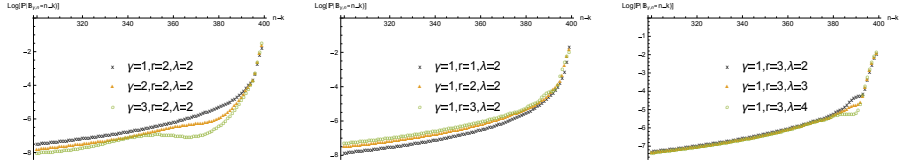

**Fig. 9** The longest block: dependency on $\gamma$ (LHS), minimum stack- (Middle) and arc-length (RHS), where $n = 400$.

## 5 The spectrum of block lengths

In the previous section we prove that a uniformly generated $\gamma$-structure almost surely contains a unique longest block. We refer to this distinguished block as the long block and any other block as short. In this section we focus on the length-distribution of short blocks.

We begin this analysis and show that with high probability we can assume, that any short block actually has finite length.

**Corollary 2** *Given any $\gamma, \epsilon > 0$, there exists an integer $t(\gamma, \epsilon)$ such that*

$$\lim_{n \to \infty} \mathbb{P}(\mathbb{B}_{\gamma, n} \geq n - t(\gamma, \epsilon)) \geq 1 - \epsilon,$$

*see Table 2.*

**Table 2** Some limit probability $\lim_{n \to \infty} \mathbb{P}(\mathbb{B}_{\gamma, n} \geq n - 100)$.

|                       | $\gamma = 1$ | $\gamma = 2$ | $\gamma = 3$ |
|-----------------------|-------------|-------------|-------------|
| $r = 2, \lambda = 2$  | 0.883       | 0.912       | 0.929       |
| $r = 3, \lambda = 4$  | 0.865       | 0.897       | 0.916       |

Next we define $g_{\gamma, k}(n, b)$ to be the number of $r$-canonical $\gamma$-structures with minimum arc-length $\lambda$, filtered by the number $b$ of blocks of length $k$. Let $\mathbf{G}_{\gamma, k}(z, u) = \sum_{n, b} g_{\gamma, k}(n, b) z^n u^b$ denote the corresponding bivariate generating function.

**Lemma 2** *The bivariate generating function of the number of $r$-canonical $\gamma$-structures with minimum arc-length $\lambda$, filtered by blocks of length $k$, is given by*

$$\mathbf{G}_{\gamma, k}(z, u) = \frac{1}{1 - \mathbf{F}_\gamma(z) - (u - 1) f_\gamma(k) z^k}.$$

The idea here is to enhance the decomposition of a $\gamma$-structure into a sequence of blocks by marking each block of length $k$. That is, we label each block of length $k$ using the term $(u - 1) f_\gamma(k) z^k$.

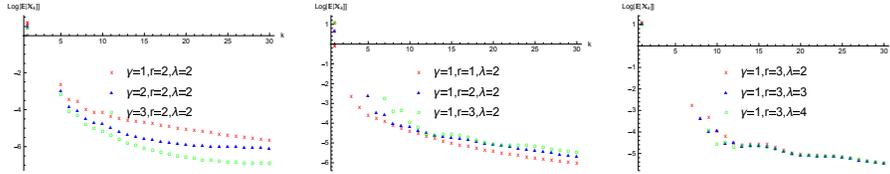

**Fig. 10** The expectation value of blocks of length $k$ for different $\gamma$, minimum stack- and arc-lengths.

Now we are in position to analyze the distribution of blocks of finite length. Let $\mathbb{X}_{\gamma,k,n}$ denote the r.v. counting the number of blocks of length $k$ in a random $\gamma$-structure over $n$ nucleotides. By construction, we have

$$\mathbb{P}(\mathbb{X}_{\gamma,k,n}=b)=\frac{g_{\gamma,k}(n,b)}{g_\gamma(n)}=\frac{[z^n u^b]\mathbf{G}_{\gamma,k}(z,u)}{[z^n]\mathbf{G}_\gamma(z)}.$$

**Theorem 4** *For fixed $\gamma$ and $k$, the distribution of the number of blocks of length $k$ in a random $\gamma$-structure of long sequence tends to a negative binomial distribution $NB(2,t)$. That is,*

$$\lim_{n\to\infty}\mathbb{P}(\mathbb{X}_{\gamma,k,n}=b)=(b+1)t^b(1-t)^2,$$

*where $\tau'=\mathbf{F}_\gamma(\rho_\gamma)$ and $t=\frac{f_\gamma(k)\rho_\gamma^k}{1-\tau'+f_\gamma(k)\rho_\gamma^k}$.*

Theorem 4 generalizes the results on secondary structures (Li and Reidys, 2018) to $\gamma$-structures, and its proof is presented in Section 8.

**Corollary 3** *For fixed $\gamma$ and $k$, the expectation of $\mathbb{X}_{\gamma,k,n}$ is asymptotically given by*

$$\lim_{n\to\infty}\mathbb{E}(\mathbb{X}_{\gamma,k,n})=\frac{2}{1-\tau'}f_\gamma(k)\rho_\gamma^k. \tag{9}$$

Fig. 10 illustrates the dependency of the expectation value of blocks of length $k$ on $\gamma$, minimum stack- and minimum arc-length. We observe that structures have fewer short blocks when structural complexity $\gamma$ increases or minimum stack-length $r$ decreases, and it is less affected by the change of minimum arc-length $\lambda$.

## 6 block types

In this section, we have a closer look at each respective type of block. By definition, a block is characterized by the irreducible shadow of its maximal component that appears in the block decomposition of the $\gamma$-structure. While the unique irreducible shadow of genus zero corresponds to blocks with a rainbow (T-type), the four irreducible shadows of genus one correspond to the following four types of blocks: H-type pseudoknots, kissing hairpins (K-type), 3-knots (L-type) and 4-knots (M-type), see Fig. 11.

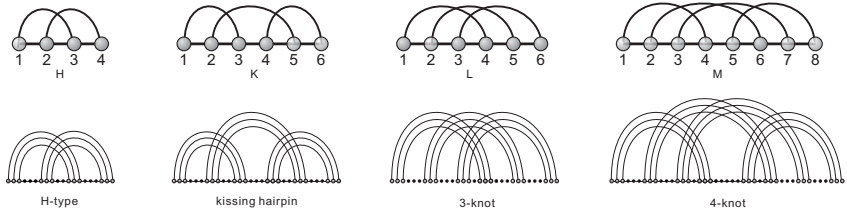

**Fig. 11** The four irreducible shadows of genus 1 and four types of blocks.

To investigate the distribution of blocks of different types in $\gamma$-structures, we consider the generating functions for blocks of type $I$. For $I \in \{T, H, K, L, M\}$, let $b_\gamma^I(n)$ denote the number of blocks of length $n$ and type $I$ in $\gamma$-structures, and $\mathbf{B}_\gamma^I(z)$ denote its corresponding generating function. We derive

**Proposition 1** *Suppose $\gamma, \lambda, r \geq 1$. The generating functions $\mathbf{B}_\gamma^I(z)$ for $I \in \{T, H, K, L, M\}$ are given by*

$$\mathbf{B}_\gamma^I(z) = h^I(z, \mathbf{G}_\gamma(z)), \tag{10}$$

*where $h^I(z, X)$ are rational functions in $z$ and $X$, given in Table 3.*

**Table 3** The rational function $h^I(z, X)$.

| T-type | H-type | K- and L-type | M-type |
|---|---|---|---|
| $\frac{z^{2r}}{1-z^2+z^{2r}}\left(X - \sum_{i=0}^{\lambda-2} z^i\right)$ | $\frac{1}{X}\left(\frac{z^{2r}X^2}{1-z^2+z^{2r}-z^{2r}X^2}\right)^2$ | $\frac{1}{X}\left(\frac{z^{2r}X^2}{1-z^2+z^{2r}-z^{2r}X^2}\right)^3$ | $\frac{1}{X}\left(\frac{z^{2r}X^2}{1-z^2+z^{2r}-z^{2r}X^2}\right)^4$ |

The proof is based on the inflation process from a particular irreducible shadow to blocks of the corresponding type, and analogous to that of Theorem 1.

Employing eq. (10), we derive the asymptotics for blocks of type $I$ from the singularity analysis of $\gamma$-structures.

**Proposition 2** *Suppose $1 \leq \gamma, \lambda, r \leq 4$ with $\lambda \leq r + 1$. The dominant singularity of $\mathbf{B}_\gamma^I(z)$ is $\rho_\gamma$, the same as that of $\mathbf{G}_\gamma(z)$. The singular expansion of $\mathbf{B}_\gamma^I(z)$ is given by*

$$\mathbf{B}_\gamma^I(z) = \tau^I + \delta^I(\rho_\gamma - z)^{\frac{1}{2}} + \theta^I(\rho_\gamma - z) + O\left((\rho_\gamma - z)^{\frac{3}{2}}\right), \qquad \text{as } z \to \rho_\gamma,$$

*where $\tau^I$, $\delta^I$ and $\theta^I$ are constants. Furthermore, the coefficients of $\mathbf{B}_\gamma^I(z)$ satisfy*

$$[z^n]\mathbf{B}_\gamma^I(z) = c^I n^{-\frac{3}{2}} \rho_\gamma^{-n} \left(1 + O(n^{-1})\right), \qquad \text{as } n \to \infty, \tag{11}$$

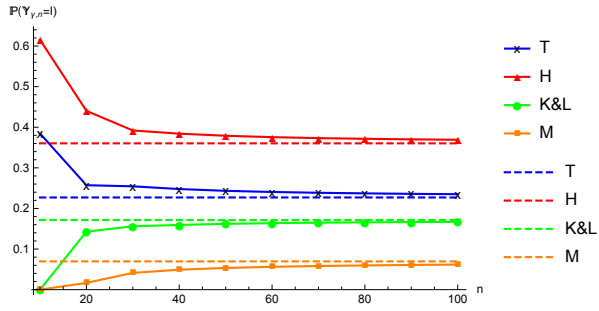

**Fig. 12** The probability of a block to be of type $I$ in 1-structures: we compare the limit probability (dashed lines) with the average probability (points) computed by exact enumeration in 1-structures. Minimum arc- and stack-length are $r = \lambda = 2$, sequence length $n = 10 \times i$, where $1 \le i \le 10$.

where $c^I$ is the positive constant $c^I = \delta^I \rho_\gamma^{\frac{1}{2}} \Gamma(-\frac{1}{2})^{-1}$.

The proof of Proposition 2 is presented in Section 8.

Now we are in position to analyze the discrete random variable, $\mathbb{Y}_{\gamma,n}$, representing the type of a block, in the space of uniformly generated blocks in $\gamma$-structures, having length $n$ with minimum stack-length $r$ and minimum arc-length $\lambda$.

**Theorem 5** *Suppose $1 \le \gamma, \lambda, r \le 4$ with $\lambda \le r + 1$. As the length of a block grows large, the probability of a block to be of type $I$ is given by*

$$\lim_{n \to \infty} \mathbb{P}(\mathbb{Y}_{\gamma,n} = I) = \delta^I \delta^{-1} \tau^2, \tag{12}$$

*where $\tau = \mathbf{G}_\gamma(\rho_\gamma)$. Furthermore, the limit probability is independent of $r$ and $\lambda$, depending only on $\gamma$ and $I$.*

The parameter independence is due to an algebraic relation satisfied by the dominant singularities of the generating functions for $\gamma$-structures, see Section 8.

Table 4 provides the limit probability of a block to be of type $I$ for $1 \le \gamma \le 3$. In Fig. 12, we contrast our limits with the exact probabilities of a block to be of type $I$ in 1-structures.

**Table 4** The limit probability $\lim_{n \to \infty} \mathbb{P}(\mathbb{Y}_{\gamma,n} = I)$ for $1 \le \gamma \le 3$.

| $\lim_{n \to \infty} \mathbb{P}(\mathbb{Y}_{\gamma,n} = I)$ | T-type | H-type | K- & L-type | M-type |
|---|---|---|---|---|
| $\gamma = 1$ | 0.227 | 0.360 | 0.171 | 0.070 |
| $\gamma = 2$ | 0.147 | 0.110 | 0.031 | 0.007 |
| $\gamma = 3$ | 0.113 | 0.057 | 0.012 | 0.002 |

Combining with Theorem 3, we have

**Corollary 4** *In the limit of long sequences, the probability of the longest block to be type $I$ is given by $\lim_{n\to\infty} \mathbb{P}(\mathbb{Y}_{\gamma,n} = I)$.*

In the following, we focus on the distribution of blocks of type $I$ and finite length. It turns out that our analysis for blocks including all types in Section 5 can straightforwardly be generalized to blocks of any type.

For fixed $k$, let $g^I_{\gamma,k}(n, b)$ denote the number of $r$-canonical $\gamma$-structures with minimum arc-length $\lambda$, filtered by the number $b$ of blocks of type $I$ and length $k$. Let $\mathbf{G}^I_{\gamma,k}(z, u) = \sum_{n,b} g^I_{\gamma,k}(n, b) z^n u^b$ denote the corresponding bivariate generating function. By labeling each block of type $I$ and length $k$ using the term $(u - 1)b^I_\gamma(k)z^k$, we obtain

$$\mathbf{G}^I_{\gamma,k}(z, u) = \frac{1}{1 - \mathbf{F}_\gamma(z) - (u - 1)b^I_\gamma(k)z^k}.$$

Next we analyze $\mathbb{X}^I_{\gamma,k,n}$, the r.v. counting the number of blocks of type $I$ and length $k$ in a random $\gamma$-structure over $n$ nucleotides. By construction, we have

$$\mathbb{P}(\mathbb{X}^I_{\gamma,k,n} = b) = \frac{g^I_{\gamma,k}(n, b)}{g_\gamma(n)} = \frac{[z^n u^b]\mathbf{G}^I_{\gamma,k}(z, u)}{[z^n]\mathbf{G}_\gamma(z)}.$$

**Theorem 6** *For fixed $\gamma$ and $k$, $\mathbb{X}^I_{\gamma,k,n}$ satisfies the discrete limit law*

$$\lim_{n\to\infty} \mathbb{P}(\mathbb{X}^I_{\gamma,k,n} = b) = (b + 1)t^b(1 - t)^2,$$

*where $\tau' = \mathbf{F}_\gamma(\rho_\gamma)$ and $t = \frac{b^I_\gamma(k)\rho_\gamma^k}{1 - \tau' + b^I_\gamma(k)\rho_\gamma^k}$. That is, the limit law of $\mathbb{X}^I_{\gamma,k,n}$ is a negative binomial distribution $NB(2, t)$ with the expectation*

$$\lim_{n\to\infty} \mathbb{E}(\mathbb{X}^I_{\gamma,k,n}) = \frac{2}{1 - \tau'}b^I_\gamma(k)\rho_\gamma^k.$$

The proof is analogous to that of Theorem 4. The key point of the proof is that the composition $\mathbf{G}^I_{\gamma,k}(z, u) = \Phi(h^I(z, u))$ belongs to the subcritical paradigm of singularity analysis (Flajolet and Sedgewick, 2009), where $h^I(z, u) = \mathbf{F}_\gamma(z) + (u - 1)b^I_\gamma(k)z^k$.

Fig. 13 illustrates the expectation value of the number of blocks of different types. We notice that the expectation for pseudoknot block type (H, K, L, M) is unimodal with respect to the length $k$, i.e., pseudoknot blocks having certain short length $k$ are more likely to be observed.

Fig. 14 shows the dependency of the expected number of blocks of length $k$ and type H on $\gamma$, minimum stack- and minimum arc-length. For blocks of any type $I$, we observe that structures have fewer blocks of type $I$ when structural complexity $\gamma$ increases or minimum stack-length $r$ decreases, and the expectation is less affected by the change of minimum arc-length $\lambda$.

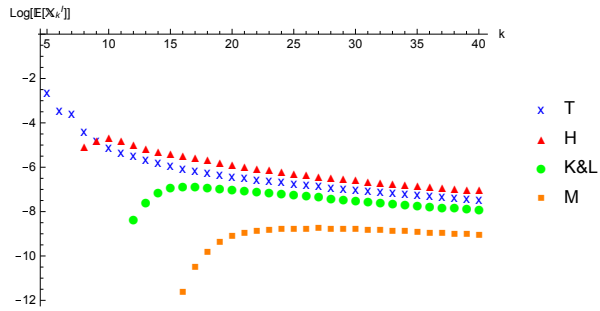

**Fig. 13** The expected number of blocks of length $k$ and different types in 1-structures with minimum arc- and stack-length $r = \lambda = 2$.

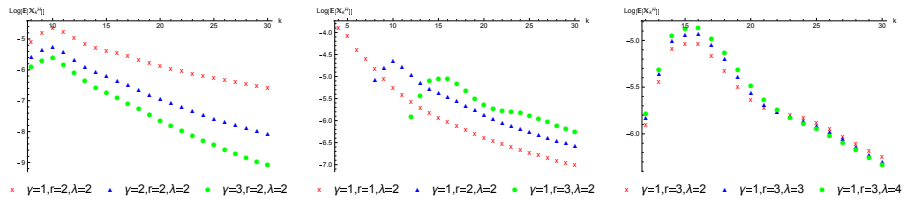

**Fig. 14** The expected number of blocks of type H and length $k$ for different $\gamma$, minimum stack- and arc-lengths.

## 7 Discussion

We have shown that the length-spectrum of blocks in random $\gamma$-structures has a gap. By Lemma 1 the longest block is a.a.s. of size $n - O(n^{1/2})$ and Corollary 2 shows that with high probability any other block has finite size. In any case, there exists a.a.s. a unique longest block. In Theorem 3 we analyze the limit distribution of the size of the unique longest block and show that it satisfies a discrete limit law. In Theorem 4 we identify the distribution of blocks of finite size $k$, in the limit of long sequences as a negative binomial. The analysis is generalized to the length-spectrum of blocks of different types. In Theorem 5, we show that the longest block is of a certain type with a fixed probability, which is independent of the choice of minimum stack- and arc-length constraint. Theorem 6 establishes the negative binomial distribution of short blocks of given type.

The analysis in Section 4 can be generalized to the lengths of the second and third longest blocks in uniformly generated $\gamma$-structures. One can show that $\mathbb{E}[\mathbb{B}_{\gamma,n}^{[2]}] = \alpha \, n^{\frac{1}{2}}\big(1 + o(1)\big)$ and $\mathbb{E}[\mathbb{B}_{\gamma,n}^{[3]}] = o(n^{\frac{1}{2}})$, where $\mathbb{B}_{\gamma,n}^{[2]}$ and $\mathbb{B}_{\gamma,n}^{[3]}$ denote the length of the second and third longest block in $\gamma$-structures, see Li and Reidys (2018).

Intuitively, one may stipulate that if a structure has more building blocks, it will likely have a larger longest block. Surprisingly, this is not true for $\gamma$-structures. As observed in Table 1, structures with greater minimum arc-length $\lambda$ have typically a larger longest block. At the same time, structures with greater $\lambda$ has fewer building blocks, i.e., the dominant singularity $\rho$ is greater. In contrast, structures with greater $\gamma$ or smaller minimum stack-length $r$ have a larger expected longest block, while the number of the building blocks in these cases grow faster.

The longest block has type $I$ with a constant probability, in the limit of long sequence (Theorem 5). On the one hand, this limit probability is invariant with respect to minimum arc- and stack-length constraint $\lambda$ and $r$, as a result of an algebraic equation for the dominant singularity. On the other hand, the limit probability decreases as structural complexity $\gamma$ increases. This is because a huge amount of new block-types flood into the structures and dilute the concentration of five specific block types we are interested in. For example, the numbers of block-types of genus 2 and 3 are 3280 and 14004032, respectively. To avoid the dilution effect and focus on the five specific block types, we present here the conditional probability for a block to be type $I$, given that the block is one of the five types, see Table 5. Remarkably, we have a higher chance to observe a block of type H than a block with a rainbow (type T) in 1-structures. For higher $\gamma$, we find more and more T-blocks.

**Table 5** The probability $\lim_{n \to \infty} \mathbb{P}(\mathbb{Y}_{\gamma,n} = I | \mathbb{Y}_{\gamma,n} \in \{\text{T, H, K, L, M}\})$ for $1 \le \gamma \le 3$.

|              | T-type | H-type | K- & L-type | M-type |
|--------------|--------|--------|-------------|--------|
| $\gamma = 1$ | 0.227  | 0.360  | 0.171       | 0.070  |
| $\gamma = 2$ | 0.450  | 0.337  | 0.095       | 0.023  |
| $\gamma = 3$ | 0.575  | 0.291  | 0.061       | 0.011  |

It is interesting to compare our results with the findings in Li and Reidys (2017) on the probability of five pseudoknot types in structures of fixed genus. Li and Reidys (2017) show that the expectation value of H-type, K-type, L-type and M-type pseudoknots, in uniformly generated structures of any genus is $O\left(n^{-1}\right)$, $O(n^{-\frac{1}{2}})$, $O(n^{-\frac{1}{2}})$ and $O(1)$, respectively. In other words, structures of genus one is dominated by M-type pseudoknots and we can hardly find H-types. In contrast, we here show that 1-structures, which generalize genus one structures by nesting and concatenating, have a much higher chance to observe H-types.

Our results are closely connected with the length of the longest arc in $\gamma$-structures. It is clear that the longest arc in a block of type $T$ has the same length as the block. According to the pigeonhole principle, the longest arc in a block of type H (or K, L, M) has at least one half (or one third) of the length of the block. Obviously, the longest arc in 1-structures is longer than the longest arc in the unique longest block. Combing with the limit probability of each block type (Theorem 5), we derive a lower bound on the expected length $\mathbb{Z}_{1,n}$ of the longest arc in 1-structures, as the sequence length grows large,

$$\mathbb{E}[\mathbb{Z}_{1,n}] \geq \sum_I d_I \lim_{n\to\infty} \mathbb{P}(\mathbb{Y}_{1,n} = I) \, \mathbb{E}[\mathbb{B}_{1,n}]$$
$$= 0.487n + o(n),$$

where $d_T = 1$, $d_H = 1/2$ and $d_K = d_L = d_M = 1/3$. This bound shows us random pseudoknot structures have the expected length of its longest arc at the same order of $n$, where we used a crude estimate $d_I$ for the length of the longest arc in each block type $I$.

Long-range arcs in RNA statistics play an important role in the context of sparsification, a particular method facilitating a speed up of the dynamic programming algorithms for the RNA folding (Wexler et al, 2007; Salari et al, 2010; Backofen et al, 2011). The theoretical analysis (Wexler et al, 2007) concludes a linear reduction time complexity based on a specific property of arcs in RNA molecules. This property is called *polymer-zeta property* and originates from studies of bonds in proteins. Polymer-zeta stipulates that long-distance base pairs have low probability. Möhl et al (2010) further apply sparsification to RNA structure prediction including pseudoknots. Our results show that in random RNA pseudoknot structures, with high probability, the longest arc has the same order of $n$. Thus the polymer-zeta property does not hold for RNA pseudoknot structure, unless one considers particular classes of natural RNA structures such as mRNA (Wexler et al, 2007).

Furthermore we have a connection between long-range blocks and the 5′-3′ distance, the length of the shortest path connecting the 5′ and 3′ ends. The finiteness of the 5′-3′ distance has first been studied in Yoffe et al (2011). Remarkably, the 5′-3′ distance of biological RNA structures is also observed to be finite, indicating that certain features of random structures can also be observed in biological structures. Han and Reidys (2012) study the $T$-blocks (rainbows) of uniformly sampled RNA secondary structures and show that the 5′-3′ distance satisfies a discrete limit law and thus is finite. Clote et al (2012) shows that the expected distance between 5′ and 3′ ends of a specific RNA sequence is finite, with respect to the Turner energy model. More importantly, the finiteness of the 5′-3′ distance and the existence of a long $T$-block both lead to the effective circularization of linear RNA, which plays an important role in many biological processes (Yoffe et al, 2011).

As for future work, we are concerned with the implications of the results of this paper for the 5′-3′ distance of $\gamma$-structures. We argue here that, using the additivity of the 5′-3′ distance with respect to our block decomposition, we can reduce the problem to the 5′-3′ distance of $\gamma$-blocks. Since $\gamma$-blocks are characterized by its irreducible shadows, it would be interesting to find out how irreducible shadows impact on the 5′-3′ distance.

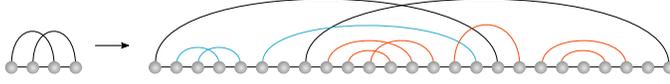

**Fig. 15** **Step I:** inflation of each arc in $\sigma$ into a sequence of induced arcs.

## 8 Proofs

*Proof (Proof of Theorem 1)* We first decompose a $\gamma$-structure into a sequence of blocks. Then we distinguish the classes of blocks into three categories: trivial blocks, 0-blocks and $\gamma$-blocks. For the latter two classes, we further decompose them by removing their maximal component, see Fig. 6. We then translate each construction into the system of functional equations in the language of symbolic enumeration.

The decomposition of $\gamma$-structures into a sequence of blocks implies eq. (1). The classification of blocks into three categories gives eq. (2).

Given a 0-block, its rainbow arc induces a maximum stack containing this arc and the nested substructure. An arc corresponds to $z^2$ and a stack of size at least $r$ corresponds to $\frac{(z^2)^r}{1-z^2}$. By construction, the nested substructure must not be a 0-block and contains at least $\lambda-1$ vertices as a result of the arc-length restriction. Since any segment with at most $\lambda-2$ vertices corresponds to the term $\sum_{i=0}^{\lambda-2} z^i$, the nested substructure gives rise to $\mathbf{G}_\gamma(z) - \mathbf{B}_0(z) - \sum_{i=0}^{\lambda-2} z^i$. Thus we arrive at eq. (3).

Given a $\gamma$-block, its maximal component can be uniquely projected into an irreducible shadow $\sigma$ of genus at most $\gamma$ having $m$ arcs. Let $\mathbf{B}_\sigma(z)$ be the generating function of blocks, having $\sigma$ as the shadow of its unique maximal component. Then we have

$$\mathbf{B}_\gamma(z) = \sum_{\sigma \in \mathcal{I}_\gamma} \mathbf{B}_\sigma(z),$$

where $\mathcal{I}_\gamma$ denotes the set of irreducible shadows of genus at most $\gamma$. We shall construct $\mathbf{B}_\gamma(z)$ in three steps using arcs, $z^2$, stacks, $\frac{(z^2)^r}{1-z^2}$, induced arcs, $\mathbf{N}(z)$, sequence of induced arcs, $\mathbf{M}(z)$, and arbitrary $\gamma$-structures, $\mathbf{G}_\gamma(z)$.

**Step I:** We inflate each arc in $\sigma$ into a sequence of induced arcs, see Fig. 15. An induced arc is an arc together with at least one nontrivial $\gamma$-structures in two intervals of both ends. Clearly, we have for a single induced arc $\mathbf{N}(z) = z^2(\mathbf{G}_\gamma(z)^2 - 1)$ and for a sequence of induced arcs $\mathbf{M}(z) = \frac{1}{1-z^2(\mathbf{G}_\gamma(z)^2-1)}$. Inflating each arc into a sequence of induced arcs produces the corresponding generating function

$$z^{2m}\mathbf{M}(z)^m = \left(\frac{z^2}{1 - z^2\left(\mathbf{G}_\gamma(z)^2 - 1\right)}\right)^m.$$

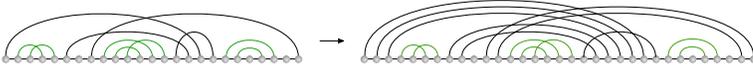

**Fig. 16 Step II:** inflation of each arc in the component with shadow $\sigma$ into stacks.

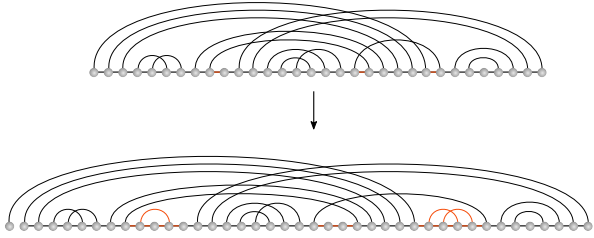

**Fig. 17 Step III:** insertion of additional $\gamma$-structures at exactly $(2m-1)$ intervals.

**Step II:** We inflate each arc in the component with shadow $\sigma$, see Fig. 16. The corresponding generating function is

$$\left(\frac{\frac{z^{2r}}{1-z^2}}{1-\frac{z^{2r}}{1-z^2}\left(\mathbf{G}_\gamma(z)^2-1\right)}\right)^m = \left(\frac{z^{2r}}{1-z^2+z^{2r}-z^{2r}\mathbf{G}_\gamma(z)^2}\right)^m$$

**Step III:** We insert additional $\gamma$-structures at exactly $(2m-1)$ intervals, see Fig. 17. Accordingly, the generating function is $\mathbf{G}_\gamma(z)^{2m-1}$.

Combining these three steps, we arrive at

$$\begin{aligned}
\mathbf{B}_\sigma(z) &= \left(\frac{z^{2r}}{1-z^2+z^{2r}-z^{2r}\mathbf{G}_\gamma(z)^2}\right)^m \mathbf{G}_\gamma(z)^{2m-1} \\
&= \mathbf{G}_\gamma(z)^{-1}\left(\frac{z^{2r}\mathbf{G}_\gamma(z)^2}{1-z^2+z^{2r}-z^{2r}\mathbf{G}_\gamma(z)^2}\right)^m.
\end{aligned}$$

Therefore

$$\begin{aligned}
\mathbf{B}_\gamma(z) &= \sum_{\sigma\in\mathcal{I}_\gamma}\mathbf{B}_\sigma(z) \\
&= \sum_{g\le\gamma}\sum_m \mathbf{i}(g,m)\,\mathbf{G}_\gamma(z)^{-1}\left(\frac{z^{2r}\mathbf{G}_\gamma(z)^2}{1-z^2+z^{2r}-z^{2r}\mathbf{G}_\gamma(z)^2}\right)^m \\
&= \mathbf{G}_\gamma(z)^{-1}\sum_{g\le\gamma}\mathbf{I}_g\Big(\frac{z^{2r}\mathbf{G}_\gamma(z)^2}{1-z^2+z^{2r}-z^{2r}\mathbf{G}_\gamma(z)^2}\Big),
\end{aligned}$$

whence eq. (4).

Note that $\sum_{g \le \gamma} \mathbf{I}_g(z)$ is a polynomial in $z$ of degree $6\gamma - 2$. Solving eqs. (1)–(4), we obtain the polynomial

$$
\begin{aligned}
Q_\gamma(z, X) =& (1 - z^2 + z^{2r}) q(z, X) \big( 1 - X + z\,X \big) \\
&+ z^{2r} q(z, X) \Big( X - \sum_{i=0}^{\lambda-2} z^i \Big) X \\
&+ q(z, X) \sum_{g \le \gamma} \mathbf{I}_g \left( \frac{z^{2r} X^2}{1 - z^2 + z^{2r} - z^{2r} X^2} \right),
\end{aligned}
$$

where $q(z, X) = (1 - z^2 + z^{2r} - z^{2r} X^2)^{6\gamma - 2}$ and $\deg(Q_\gamma(z, X)) = 12\gamma - 2$. Then $Q_\gamma(z, \mathbf{G}_\gamma(z)) = 0$, completing the proof.

*Proof (Sketch proof of Lemma 1)* (Li and Reidys, 2018) The proof is based on three claims. First, we provide an asymptotic formula for the probability $\mathbb{P}(\mathbb{B}_{\gamma,n} = n - k)$ as $k, n \to \infty$ and $k < \frac{n}{2}$. This will motivate us to partition the expectation $\mathbb{E}[\mathbb{B}_{\gamma,n}]$ into two sums, i.e.,

$$
\begin{aligned}
\mathbb{E}[\mathbb{B}_{\gamma,n}] &= \sum_{k=0}^{n-1} (n - k) \mathbb{P}(\mathbb{B}_{\gamma,n} = n - k) \\
&= n - \sum_{0 \le k < \frac{n}{2}} k \mathbb{P}(\mathbb{B}_{\gamma,n} = n - k) - \sum_{\frac{n}{2} \le k \le n-1} k \mathbb{P}(\mathbb{B}_{\gamma,n} = n - k).
\end{aligned}
$$

We proceed by estimating $\sum_{0 \le k < \frac{n}{2}} k \mathbb{P}(\mathbb{B}_{\gamma,n} = n-k)$ and $\sum_{\frac{n}{2} \le k \le n-1} k \mathbb{P}(\mathbb{B}_{\gamma,n} = n - k)$.

**Claim 1:** As $k, n \to \infty$ with $k < \frac{n}{2}$, we have

$$
\mathbb{P}(\mathbb{B}_{\gamma,n} = n - k) = 2c\,\tau^{-1} \Big( 1 - \frac{k}{n} \Big)^{-\frac{3}{2}} k^{-\frac{3}{2}} \big( 1 + O(k^{-1}) \big) \big( 1 + O(n^{-1}) \big), \quad (13)
$$

where $c$ and $\tau$ are given by Theorem 2.

To derive an expression for $\mathbb{P}(\mathbb{B}_{\gamma,n} = n - k)$ in terms of coefficients of known generating functions, we define $\mathbf{F}_{\le m}(z) = \sum_{1 \le i \le m} f(i) z^i$ to be the truncated series for blocks of length at most $m$. Then the generating function for $\gamma$-structures consisting of blocks of length at most $m$ is given by $\mathbf{G}_{\le m}(z) = \frac{1}{1 - \mathbf{F}_{\le m}(z)}$. Therefore, by construction,

$$
\mathbb{P}(\mathbb{B}_{\gamma,n} = n - k) = \frac{[z^n](\mathbf{G}_{\le n-k}(z) - \mathbf{G}_{\le n-k-1}(z))}{[z^n]\mathbf{G}_\gamma(z)}.
$$

To proceed we consider the Taylor expansion of $\mathbf{G}_{\le n-k}(z)$

$$
\mathbf{G}_{\le n-k}(z) = \Phi(\mathbf{F}_{\le n-k}(z)) = \sum_{i \ge 0} \frac{(-1)^i \Phi^{(i)}(\mathbf{F}_\gamma(z))}{i!} \left( \mathbf{F}_\gamma(z) - \mathbf{F}_{\le n-k}(z) \right)^i,
$$
$$(14)$$

where $\Phi(z) = \frac{1}{1-z}$. Eq. (14) can also be proved combinatorially using the inclusion-exclusion principle, by observing $\frac{\Phi^{(i)}(\mathbf{F}_\gamma(z))}{i!} \left( \mathbf{F}_\gamma(z) - \mathbf{F}_{\le n-k}(z) \right)^i$ counts structures with $i$ marked blocks of length greater than $n - k$. Note that $[z^n]\frac{\Phi^{(i)}(\mathbf{F}_\gamma(z))}{i!}\,(\mathbf{F}_\gamma(z) - \mathbf{F}_{\leq n-k}(z))^i = 0$ for $i \geq 2$, since $k < \frac{n}{2}$ and any structure of length $n$ has at most one block of length greater than $\frac{n}{2}$. By taking the coefficient of $z^n$ in eq. (14), we obtain

$$[z^n]\mathbf{G}_{\leq n-k}(z) = [z^n]\Big(\Phi(\mathbf{F}_\gamma(z)) + \Phi'(\mathbf{F}_\gamma(z))(\mathbf{F}_{\leq n-k}(z) - \mathbf{F}_\gamma(z))\Big). \qquad (15)$$

Similarly, eq. (15) holds for $[z^n]\mathbf{G}_{\leq n-k-1}(z)$. Therefore, we arrive at

$$
\begin{aligned}
\mathbb{P}(\mathbb{B}_{\gamma,n} = n - k) &= \frac{[z^n](\Phi'(\mathbf{F}_\gamma(z))(\mathbf{F}_{\leq n-k}(z) - \mathbf{F}_{\leq n-k-1}(z)))}{[z^n]\mathbf{G}_\gamma(z)} \\
&= \frac{[z^k]\Phi'(\mathbf{F}_\gamma(z))\,[z^{n-k}]\mathbf{F}_\gamma(z)}{[z^n]\mathbf{G}_\gamma(z)}.
\end{aligned}
\qquad (16)
$$

Combining with the singularity analysis of $\mathbf{F}_\gamma(z)$ and $\mathbf{G}_\gamma(z)$ in Theorem 2 and Corollary 1, we derive eq. (13).

**Claim 2:**

$$\sum_{0 \leq k < \frac{n}{2}} k\mathbb{P}(\mathbb{B}_{\gamma,n} = n - k) = \alpha\, n^{\frac{1}{2}}\big(1 + o(1)\big), \qquad \text{as } n \to \infty. \qquad (17)$$

In view of the fact that $\sum_{0 \leq k \leq n^{\frac{1}{8}}} k\mathbb{P}(\mathbb{B}_{\gamma,n} = n - k) = O(n^{\frac{1}{8}} \cdot n^{\frac{1}{8}}) = o(n^{\frac{1}{2}})$, this motivates to split the summation of eq. (17) and to consider the term $\sum_{n^{\frac{1}{8}} \leq k < \frac{n}{2}} k\mathbb{P}(\mathbb{B}_{\gamma,n} = n - k)$ separately. Employing eq. (13), we estimate $\sum_{n^{\frac{1}{8}} \leq k < \frac{n}{2}} k\mathbb{P}(\mathbb{B}_{\gamma,n} = n - k)$ via the integral $\int_0^{\frac{1}{2}}(1 - x)^{-\frac{3}{2}}x^{-\frac{1}{2}}\mathrm{d}x$. This leads to eq. (17), where $\alpha = 2\,c\,\tau^{-1}\int_0^{\frac{1}{2}}(1 - x)^{-\frac{3}{2}}x^{-\frac{1}{2}}\mathrm{d}x = 4\,c\,\tau^{-1}$.

**Claim 3:**

$$\sum_{\frac{n}{2} \leq k \leq n-1} k\mathbb{P}(\mathbb{B}_{\gamma,n} = n - k) = o(n^{\frac{1}{2}}), \qquad \text{as } n \to \infty. \qquad (18)$$

We compute

$$
\begin{aligned}
&\sum_{\frac{n}{2} \leq k \leq n-1} k\mathbb{P}(\mathbb{B}_{\gamma,n} = n - k) \\
&\leq n \sum_{\frac{n}{2} \leq k \leq n-1} \mathbb{P}(\mathbb{B}_{\gamma,n} = n - k) \\
&= n\Big(1 - \sum_{0 \leq k < n^{\frac{2}{8}}} \mathbb{P}(\mathbb{B}_{\gamma,n} = n - k) - \sum_{n^{\frac{2}{8}} \leq k < \frac{n}{2}} \mathbb{P}(\mathbb{B}_{\gamma,n} = n - k)\Big).
\end{aligned}
\qquad (19)
$$

Then we proceed by computing the two sums $\sum_{0 \leq k < n^{\frac{2}{8}}} \mathbb{P}(\mathbb{B}_{\gamma,n} = n - k)$ and $\sum_{n^{\frac{2}{8}} \leq k < \frac{n}{2}} \mathbb{P}(\mathbb{B}_{\gamma,n} = n - k)$.

For $0 \le k < n^{\frac{2}{5}}$, we derive from eq. (16) that

$$\mathbb{P}(\mathbb{B}_{\gamma,n} = n - k) = \frac{[z^k]\Phi'(\mathbf{F}_\gamma(z))\,[z^{n-k}]\mathbf{F}_\gamma(z)}{[z^n]\mathbf{G}_\gamma(z)} = \tau^{-2}\,b_k\,\rho_\gamma^k(1 + o(n^{-\frac{1}{2}})), \quad (20)$$

where $b_k = [z^k]\Phi'(\mathbf{F}_\gamma(z))$. Thus, we obtain

$$\sum_{0 \le k < n^{\frac{2}{5}}} \mathbb{P}(\mathbb{B}_{\gamma,n} = n - k)$$

$$= \Big(1 - \tau^{-2} \sum_{k \ge n^{\frac{2}{5}}} b_k\,\rho_\gamma^k\Big)\big(1 + o(n^{-\frac{1}{2}})\big)$$

$$= \Big(1 - 2c\,\tau^{-1} \sum_{k \ge n^{\frac{2}{5}}} k^{-\frac{3}{2}}\big(1 + O(k^{-1})\big)\Big)\big(1 + o(n^{-\frac{1}{2}})\big) \qquad (21)$$

$$= \Big(1 - 4c\,\tau^{-1}\,n^{-\frac{1}{5}}\big(1 + O(n^{-\frac{2}{5}})\big)\Big)\big(1 + o(n^{-\frac{1}{2}})\big) \qquad (22)$$

$$= 1 - \alpha n^{-\frac{1}{5}} + o(n^{-\frac{1}{2}}), \qquad \text{as } n \to \infty.$$

Eq. (21) employs the asymptotics for $b_k = [z^k]\Phi'(\mathbf{F}_\gamma(z))$. Eq. (22) is derived from the asymptotic expansion of the Hurwitz-Zeta function $\zeta(s,n) = \frac{n^{1-s}}{s-1}\big(1 + O(n^{-1})\big)$, where $\zeta(s,n) = \sum_{i=0}^{\infty}(n+i)^{-s}$.

As for the second sum, we have

$$\sum_{n^{\frac{2}{5}} \le k < \frac{n}{2}} \mathbb{P}(\mathbb{B}_{\gamma,n} = n - k)$$

$$= 2c\,\tau^{-1} \sum_{n^{\frac{2}{5}} \le k < \frac{n}{2}} \Big(1 - \frac{k}{n}\Big)^{-\frac{3}{2}} k^{-\frac{3}{2}}\big(1 + O(k^{-1})\big)\big(1 + O(n^{-1})\big) \qquad (23)$$

$$= \frac{\alpha}{2} \cdot 2n^{-\frac{1}{5}}(1 + O(n^{-\frac{2}{5}}))\big(1 + O(n^{-1})\big) \qquad (24)$$

$$= \alpha n^{-\frac{1}{5}} + o(n^{-\frac{1}{2}}), \qquad \text{as } n \to \infty.$$

Eq. (23) follows from eq. (13). In eq. (24), the summation is approximated by the Euler-Maclaurin summation formula (see, for example, Graham et al (1994)).

Accordingly, eq. (24) is established. Combining the two sums, we derive eq. (18).

Now we are in position to compute

$$\mathbb{E}[\mathbb{B}_{\gamma,n}] = n - \sum_{0 \le k < \frac{n}{2}} k\mathbb{P}(\mathbb{B}_{\gamma,n} = n - k) - \sum_{\frac{n}{2} \le k \le n-1} k\mathbb{P}(\mathbb{B}_{\gamma,n} = n - k)$$

$$= n - \alpha\,n^{\frac{1}{2}}\big(1 + o(1)\big).$$

As for the variance $\mathbb{V}[\mathbb{B}_{\gamma,n}] = \mathbb{E}[\mathbb{B}_{\gamma,n}^2] - \mathbb{E}[\mathbb{B}_{\gamma,n}]^2$, the key step is to show $\sum_{k=0}^{n-1} k^2 \mathbb{P}(\mathbb{B}_{\gamma,n} = n - k) = \beta n^{\frac{3}{2}}(1 + o(1))$, following the same line of arguments as in Claims 2 and 3.

*Proof (Proof of Theorem 3)* (Li and Reidys, 2018) According to Lemma 1, $\mathbb{B}_{\gamma,n}$ is concentrated at $n - \alpha\, n^{\frac{1}{2}}$ with a variance of $O(n^{\frac{3}{2}})$. Chebyshev's inequality then guarantees

$$\mathbb{P}\Big(\mathbb{E}[\mathbb{B}_{\gamma,n}] - \mathbb{B}_{\gamma,n} \geq a\Big) \leq \frac{\mathbb{V}[\mathbb{B}_{\gamma,n}]}{a^2}.$$

Accordingly, for $a = \Omega(n^t)$ with $t > \frac{3}{4}$, the right hand-side tends to zero as $n$ tends to infinity, whence eq. (7). To establish eq. (8) we inspect that the proof of eq. (20) in Lemma 1 holds for $k = o(n)$ and eq. (8) follows. In summary, eq. (7) implies that a.a.s. we may assume $k = o(n)$ in which case eq. (8) guarantees that $n - \mathbb{B}_{\gamma,n} = k$ satisfies a discrete limit law.

*Proof (Proof of Theorem 4)* (Li and Reidys, 2018) In view of Lemma 2, $\mathbf{G}_{\gamma,k}(z,u)$ can be expressed as $\mathbf{G}_{\gamma,k}(z,u) = \Phi(h(z,u))$, where $\Phi(z) = \frac{1}{1-z}$ and $h(z,u) = \mathbf{F}_{\gamma}(z) + (u-1)f_{\gamma}(k)z^k$. Since $h(z,u)$ have nonnegative coefficients and $h(0,0) = 0$, the composition $\Phi(h(z,u))$ is well-defined.

We verify that $\mathbf{G}_{\gamma,k}(z,u)$ has the same dominant singularity $\rho_{\gamma}$ as $\mathbf{F}_{\gamma}(z)$, by checking that there exists a neighborhood $U$ of 1 such that $h(\rho_{\gamma}, u) < 1$ for all $u$ in $U$. As a result, the composition $\mathbf{G}_{\gamma,k}(z,u) = \Phi(h(z,u))$ belongs to the subcritical case of singularity analysis (Flajolet and Sedgewick, 2009). The singular expansion of $h(z,u)$ at $\rho_{\gamma}$ is derived from that of $\mathbf{F}_{\gamma}(z)$ in Corollary 1

$$h(z,u) = \tau' + (u-1)f_{\gamma}(k)\rho_{\gamma}^k + \delta'\big(\rho_{\gamma} - z\big)^{\frac{1}{2}}(1 + o(1)).$$

Combining this with the regular expansion of $\Phi(z)$ at $\tau_1 = \tau' + (u-1)f_{\gamma}(k)\rho_{\gamma}^k$

$$\Phi(z) = \Phi(\tau_1) + \Phi'(\tau_1)(z - \tau_1)(1 + o(1)),$$

we derive the singular expansion of $\mathbf{G}_{\gamma,k}(z,u)$ at $\rho_{\gamma}$

$$\mathbf{G}_{\gamma,k}(z,u) = \Phi(\tau_1) + \Phi'(\tau_1)\delta'\big(\rho_{\gamma} - z\big)^{\frac{1}{2}}(1 + o(1)).$$

By the transfer theorem (Flajolet and Sedgewick, 2009), we obtain

$$[z^n]\mathbf{G}_{\gamma,k}(z,u) = \Phi'(\tau_1)\delta' c_k\, n^{-\frac{3}{2}}\rho^n(1 + o(1)).$$

Now we are in position to compute

$$p_k(u) = \lim_{n \to \infty} \sum_b \mathbb{P}(\mathbb{X}_{\gamma,k,n} = b)u^b$$

$$= \lim_{n \to \infty} \frac{[z^n]\mathbf{G}_{\gamma,k}(z,u)}{[z^n]\mathbf{G}_{\gamma}(z)}$$

$$= \lim_{n \to \infty} \frac{\Phi'(\tau_1)}{\Phi'(\tau')}$$

$$= \Big(\frac{1-t}{1-tu}\Big)^2,$$

where $t = \frac{f_{\gamma}(k)\rho_{\gamma}^k}{1-\tau'+f_{\gamma}(k)\rho_{\gamma}^k}$, completing the proof.

*Proof (Proof of Proposition 2)* We verify that $\mathbf{B}_\gamma^I(z)$ has the same dominant singularity $\rho_\gamma$ as $\mathbf{G}_\gamma(z)$, by checking that $1 - z^2 + z^{2r} \neq 0$ (T-type) and $1 - z^2 + z^{2r} - z^{2r}\mathbf{G}_\gamma(z)^2 \neq 0$ (H-, K-, L-, M- types) for all $|z| \leq \rho_\gamma$. As a result, the composition $\mathbf{B}_\gamma^I(z) = h^I(z, \mathbf{G}_\gamma(z))$ belongs to the subcritical paradigm of singularity analysis (Flajolet and Sedgewick, 2009). Combining the Taylor expansion of $h^I(z, X)$ with the singular expansion of $\mathbf{G}_\gamma(z)$ in Theorem 2, we obtain the singular expansion of $\mathbf{B}_\gamma^I(z)$ at $\rho_\gamma$. The asymptotics for the coefficients then follows from the transfer theorem (Flajolet and Sedgewick, 2009).

*Proof (Proof of Theorem 5)* We derive eq. (12), employing eqs. (6) and (11) and the definition $\mathbb{P}(\mathbb{Y}_{\gamma,n} = I) = \frac{[z^n]\mathbf{B}_\gamma^I(z)}{[z^n]\mathbf{F}_\gamma(z)}$.

It remains to show that $\lim_{n\to\infty}\mathbb{P}(\mathbb{Y}_{\gamma,n} = I) = \delta^I\delta^{-1}\tau^2$ is independent of $r$ and $\lambda$. By Proposition 2, we employ the subcritical paradigm to compute $\delta^I = \delta\left.\frac{\partial}{\partial X}h^I(z,X)\right|_{z=\rho_\gamma, X=\tau}$. As a result, we express the limit probability using $\rho_\gamma$ and $\tau = \mathbf{G}_\gamma(\rho_\gamma)$, see Table 6. Furthermore, we observe that it can be rewritten in terms of $\eta_\gamma$, where $\eta_\gamma = \frac{\rho_\gamma^{2r}\tau^2}{1-\rho_\gamma^2+\rho_\gamma^{2r}}$. Therefore it suffices to show $\eta_\gamma$ is independent of $r$ and $\lambda$.

**Table 6** The limit probability $\lim_{n\to\infty}\mathbb{P}(\mathbb{Y}_{\gamma,n} = I)$.

| T-type | H-type | K- and L-type | M-type |
|---|---|---|---|
| $\frac{\rho_\gamma^{2r}\tau^2}{1-\rho_\gamma^2+\rho_\gamma^{2r}}$ | $\frac{\rho_\gamma^{4r}\tau^4(3-3\rho_\gamma^2+3\rho_\gamma^{2r}+\rho_\gamma^{2r}\tau^2)}{(1-\rho_\gamma^2+\rho_\gamma^{2r}-\rho_\gamma^{2r}\tau^2)^3}$ | $\frac{\rho_\gamma^{6r}\tau^6(5-5\rho_\gamma^2+5\rho_\gamma^{2r}+\rho_\gamma^{2r}\tau^2)}{(1-\rho_\gamma^2+\rho_\gamma^{2r}-\rho_\gamma^{2r}\tau^2)^4}$ | $\frac{\rho_\gamma^{8r}\tau^8(7-7\rho_\gamma^2+7\rho_\gamma^{2r}+\rho_\gamma^{2r}\tau^2)}{(1-\rho_\gamma^2+\rho_\gamma^{2r}-\rho_\gamma^{2r}\tau^2)^5}$ |
| $\eta_\gamma$ | $\frac{\eta_\gamma^2(3+\eta_\gamma)}{(1-\eta_\gamma)^3}$ | $\frac{\eta_\gamma^3(5+\eta_\gamma)}{(1-\eta_\gamma)^4}$ | $\frac{\eta_\gamma^4(7+\eta_\gamma)}{(1-\eta_\gamma)^5}$ |

To this end, we employ a form of the generating function $\mathbf{G}_\gamma(z)$ derived in Han et al (2014), that is,

$$\mathbf{G}_\gamma(z) = \frac{1}{v(z)}\,\mathbf{H}_\gamma\left(\frac{u(z)}{v(z)^2}\right), \qquad (25)$$

where $u(z) = \frac{z^{2r}}{1-z^2+z^{2r}}$, $v(z) = 1 - z + u(z)\sum_{i=0}^{\lambda-2}z^i$ and $\mathbf{H}_\gamma(z)$ denotes the generating function of $\gamma$-matchings, i.e., $\gamma$-structures without isolated vertices and minimum stack- and arc-length constraints. By definition, both $\mathbf{H}_\gamma(z)$ and its dominant singularity $\rho_\mathbf{H}$ are independent of $r$ and $\lambda$. Han et al (2014) prove that the dominant singularity of $\mathbf{G}_\gamma(z)$ satisfies

$$\frac{u(\rho_\gamma)}{v(\rho_\gamma)^2} = \rho_\mathbf{H}.$$

Setting $z = \rho_\gamma$ in eq. (25), we have $\tau = \mathbf{G}_\gamma(\rho_\gamma) = \frac{1}{v(\rho_\gamma)}\,\mathbf{H}_\gamma(\rho_\mathbf{H})$. Therefore we derive

$$\eta_\gamma = u(\rho_\gamma)\tau^2 = \frac{u(\rho_\gamma)\mathbf{H}_\gamma(\rho_\mathbf{H})^2}{v(\rho_\gamma)^2} = \rho_\mathbf{H}\mathbf{H}_\gamma(\rho_\mathbf{H})^2,$$

implying that $\eta_\gamma$ is independent of $r$ and $\lambda$.

**ACKNOWLEDGMENTS** We want to thank Christopher Barrett for stimulating discussions and the staff of the Biocomplexity Institute of Virginia Tech for their great support. The third author is a Thermo Fisher Scientific Fellow in Advanced Systems for Information Biology and acknowledges their support of this work.

## AUTHOR DISCLOSURE STATEMENT

The authors declare that no competing financial interests exist.